\documentclass[12pt]{article}
\usepackage{latexsym}
\usepackage{amsfonts}
\usepackage{amsmath}
\usepackage{amssymb}
\usepackage{amsthm}
\newtheorem{theo}{Theorem}[subsection]

\newtheorem{lem}[theo]{Lemma}

\newtheorem{examp}[theo]{Example}

\date{}
\begin{document}
\newcommand{\supp}{{\rm supp}}
\def\ab{{\rm ab}}
\def\and{{\rm and}}
\def\cd{{\rm cd}}
\def\alg{{\rm alg}}
\def\Aut{{\rm Aut}}
\def\gcd{{\rm gcd}}
\def\GL{{\rm GL}}
\def\Hom{{\rm Hom}}
\def\id{{\rm id}}
\def\Im{{\rm Im}}
\def\Inn{{\rm Inn}}
\def\irr{{\rm irr}}
\def\Ker{{\rm Ker}}
\def\mod{\;\mathop{\rm mod}\hskip0.2em\relax}
\def\ord{{\rm ord}}
\def\Out{{\rm Out}}
\def\PGL{{\rm PGL}}
\def\PSL{{\rm PSL}}
\def\pr{{\rm pr}}
\def\rank{{\rm rank}}
\def\ring{{\rm ring}}
\def\res{{\rm res}}
\def\Res{{\rm Res}}
\def\SL{{\rm SL}}
\def\Spec{{\rm Spec}}
\def\Subg{{\rm Subg}}
\newcommand{\calA}{{\mathcal A}}
\newcommand{\calB}{{\mathcal B}}
\newcommand{\calC}{{\mathcal C}}
\newcommand{\calD}{{\mathcal D}}
\newcommand{\calE}{{\mathcal E}}
\newcommand{\calF}{{\mathcal F}}
\newcommand{\calG}{{\mathcal G}}
\newcommand{\calH}{{\mathcal H}}
\newcommand{\calI}{{\mathcal I}}
\newcommand{\calJ}{{\mathcal J}}
\newcommand{\calK}{{\mathcal K}}
\newcommand{\calL}{{\mathcal L}}
\newcommand{\calM}{{\mathcal M}}
\newcommand{\calN}{{\mathcal N}}
\newcommand{\calO}{{\mathcal O}}
\newcommand{\calP}{{\mathcal P}}
\newcommand{\calQ}{{\mathcal Q}}
\newcommand{\calR}{{\mathcal R}}
\newcommand{\calS}{{\mathcal S}}
\newcommand{\calT}{{\mathcal T}}
\newcommand{\calU}{{\mathcal U}}
\newcommand{\calV}{{\mathcal V}}
\newcommand{\calW}{{\mathcal W}}
\newcommand{\calX}{{\mathcal X}}
\newcommand{\calY}{{\mathcal Y}}
\newcommand{\calZ}{{\mathcal Z}}
\newcommand{\agal}{{\tilde a}}       \newcommand{\Agal}{{\tilde A}}
\newcommand{\bgal}{{\tilde b}}       \newcommand{\Bgal}{{\tilde B}}
\newcommand{\cgal}{{\tilde c}}       \newcommand{\Cgal}{{\tilde C}}
\newcommand{\dgal}{{\tilde d}}       \newcommand{\Dgal}{{\tilde D}}
\newcommand{\egal}{{\tilde e}}       \newcommand{\Egal}{{\tilde E}}
\newcommand{\fgal}{{\tilde f}}       \newcommand{\Fgal}{{\tilde F}}
\newcommand{\ggal}{{\tilde g}}       \newcommand{\Ggal}{{\tilde G}}
\newcommand{\hgal}{{\tilde h}}       \newcommand{\Hgal}{{\tilde H}}
\newcommand{\igal}{{\tilde i}}       \newcommand{\Igal}{{\tilde I}}
\newcommand{\jgal}{{\tilde j}}       \newcommand{\Jgal}{{\tilde J}}
\newcommand{\kgal}{{\tilde k}}       \newcommand{\Kgal}{{\tilde K}}
\newcommand{\lgal}{{\tilde l}}       \newcommand{\Lgal}{{\tilde L}}
\newcommand{\mgal}{{\tilde m}}       \newcommand{\Mgal}{{\tilde M}}
\newcommand{\ngal}{{\tilde n}}       \newcommand{\Ngal}{{\tilde N}}
\newcommand{\ogal}{{\tilde o}}       \newcommand{\Ogal}{{\tilde O}}
\newcommand{\pgal}{{\tilde p}}       \newcommand{\Pgal}{{\tilde P}}
\newcommand{\qgal}{{\tilde q}}       \newcommand{\Qgal}{{\tilde Q}}
\newcommand{\rgal}{{\tilde r}}       \newcommand{\Rgal}{{\tilde R}}
\newcommand{\sgal}{{\tilde s}}       \newcommand{\Sgal}{{\tilde S}}
\newcommand{\tgal}{{\tilde t}}       \newcommand{\Tgal}{{\tilde T}}
\newcommand{\ugal}{{\tilde u}}       \newcommand{\Ugal}{{\tilde U}}
\newcommand{\vgal}{{\tilde v}}       \newcommand{\Vgal}{{\tilde V}}
\newcommand{\wgal}{{\tilde w}}       \newcommand{\Wgal}{{\tilde W}}
\newcommand{\xgal}{{\tilde x}}       \newcommand{\Xgal}{{\tilde X}}
\newcommand{\ygal}{{\tilde y}}       \newcommand{\Ygal}{{\tilde Y}}
\newcommand{\zgal}{{\tilde z}}       \newcommand{\Zgal}{{\tilde Z}}
\newcommand{\gag}[1]{\bar{#1}}
\def\agag{{\bar a}}     \def\Agag{{\gag{A}}}
\def\bgag{{\bar b}}     \def\Bgag{{\gag{B}}}
\def\cgag{{\bar c}}     \def\Cgag{{\gag{C}}}
\def\dgag{{\bar d}}     \def\Dgag{{\gag{D}}}
\def\egag{{\bar e}}     \def\Egag{{\gag{E}}}
\def\fgag{{\bar f}}     \def\Fgag{{\gag{F}}}
\def\ggag{{\bar g}}     \def\Ggag{{\gag{G}}}
\def\hgag{{\bar h}}     \def\Hgag{{\gag{H}}}
\def\igag{{\bar i}}     \def\Igag{{\gag{I}}}
\def\jgag{{\bar \jmath}}    \def\Jgag{{\bar{J}}}
\def\kgag{{\bar k}}     \def\Kgag{{\gag{K}}}
\def\lgag{{\bar l}}     \def\Lgag{{\gag{L}}}
\def\mgag{{\bar m}}     \def\Mgag{{\gag{M}}}
\def\ngag{{\bar n}}     \def\Ngag{{\gag{N}}}
\def\ogag{{\bar o}}     \def\Ogag{{\gag{O}}}
\def\pgag{{\bar p}}     \def\Pgag{{\gag{P}}}
\def\qgag{{\bar q}}     \def\Qgag{{\gag{Q}}}
\def\rgag{{\bar r}}     \def\Rgag{{\gag{R}}}
\def\sgag{{\bar s}}     \def\Sgag{{\gag{S}}}
\def\tgag{{\bar t}}     \def\Tgag{{\gag{T}}}
\def\ugag{{\bar u}}     \def\Ugag{{\gag{U}}}
\def\vgag{{\bar v}}     \def\Vgag{{\gag{V}}}
\def\wgag{{\bar w}}     \def\Wgag{{\gag{W}}}
\def\xgag{{\bar x}}     \def\Xgag{{\gag{X}}}
\def\ygag{{\bar y}}     \def\Ygag{{\gag{Y}}}
\def\zgag{{\bar z}}     \def\Zgag{{\gag{Z}}}
\def\ahat{{\hat a}}     \def\Ahat{{\hat A}}
\def\bhat{{\hat b}}     \def\Bhat{{\hat B}}
\def\chat{{\hat c}}     \def\Chat{{\hat C}}
\def\dhat{{\hat d}}     \def\Dhat{{\hat D}}
\def\ehat{{\hat e}}     \def\Ehat{{\hat E}}
\def\fhat{{\hat f}}     \def\Fhat{{\hat F}}
\def\ghat{{\hat g}}     \def\Ghat{{\hat G}}
\def\hhat{{\hat h}}     \def\Hhat{{\hat H}}
\def\ihat{{\hat i}}     \def\Ihat{{\hat I}}
\def\jhat{{\hat j}}     \def\Jhat{{\hat J}}
\def\khat{{\hat k}}     \def\Khat{{\hat K}}
\def\lhat{{\hat l}}     \def\Lhat{{\hat L}}
\def\mhat{{\hat m}}     \def\Mhat{{\hat M}}
\def\nhat{{\hat n}}     \def\Nhat{{\hat N}}
\def\ohat{{\hat o}}     \def\Ohat{{\hat O}}
\def\phat{{\hat p}}     \def\Phat{{\hat P}}
\def\qhat{{\hat q}}     \def\Qhat{{\hat Q}}
\def\rhat{{\hat r}}     \def\Rhat{{\hat R}}
\def\shat{{\hat s}}     \def\Shat{{\hat S}}
\def\that{{\hat t}}     \def\That{{\hat T}}
\def\uhat{{\hat u}}     \def\Uhat{{\hat U}}
\def\vhat{{\hat v}}     \def\Vhat{{\hat V}}
\def\what{{\hat w}}     \def\What{{\hat W}}
\def\xhat{{\hat x}}     \def\Xhat{{\hat X}}
\def\yhat{{\hat y}}     \def\Yhat{{\hat Y}}
\def\zhat{{\hat z}}     \def\Zhat{{\hat Z}}
\newcommand{\bbA}{\mathbb{A}}
\newcommand{\bbB}{\mathbb{B}}
\newcommand{\bbC}{\mathbb{C}}
\newcommand{\bbD}{\mathbb{D}}
\newcommand{\bbE}{\mathbb{E}}
\newcommand{\bbF}{\mathbb{F}}
\newcommand{\bbG}{\mathbb{G}}
\newcommand{\bbH}{\mathbb{H}}
\newcommand{\bbI}{\mathbb{I}}
\newcommand{\bbJ}{\mathbb{J}}
\newcommand{\bbK}{\mathbb{K}}
\newcommand{\bbL}{\mathbb{L}}
\newcommand{\bbM}{\mathbb{M}}
\newcommand{\bbN}{\mathbb{N}}
\newcommand{\bbO}{\mathbb{O}}
\newcommand{\bbP}{\mathbb{P}}
\newcommand{\bbQ}{\mathbb{Q}}
\newcommand{\bbR}{\mathbb{R}}
\newcommand{\bbS}{\mathbb{S}}
\newcommand{\bbT}{\mathbb{T}}
\newcommand{\bbU}{\mathbb{U}}
\newcommand{\bbV}{\mathbb{V}}
\newcommand{\bbW}{\mathbb{W}}
\newcommand{\bbX}{\mathbb{X}}
\newcommand{\bbY}{\mathbb{Y}}
\newcommand{\bbZ}{\mathbb{Z}}
\newcommand{\bfa}{\mbox{\boldmath$a$}}
\newcommand{\bfb}{\mbox{\boldmath$b$}}
\newcommand{\bfc}{\mbox{\boldmath$c$}}
\newcommand{\bfd}{\mbox{\boldmath$d$}}
\newcommand{\bfe}{\mbox{\boldmath$e$}}
\newcommand{\bff}{\mbox{\boldmath$f$}}
\newcommand{\bfg}{\mbox{\boldmath$g$}}
\newcommand{\bfh}{\mbox{\boldmath$h$}}
\newcommand{\bfi}{\mbox{\boldmath$i$}}
\newcommand{\bfj}{\mbox{\boldmath$j$}}
\newcommand{\bfk}{\mbox{\boldmath$k$}}
\newcommand{\bfl}{\mbox{\boldmath$l$}}
\newcommand{\bfm}{\mbox{\boldmath$m$}}
\newcommand{\bfn}{\mbox{\boldmath$n$}}
\newcommand{\bfo}{\mbox{\boldmath$o$}}
\newcommand{\bfp}{\mbox{\boldmath$p$}}
\newcommand{\bfq}{\mbox{\boldmath$q$}}
\newcommand{\bfr}{\mbox{\boldmath$r$}}
\newcommand{\bfs}{\mbox{\boldmath$s$}}
\newcommand{\bft}{\mbox{\boldmath$t$}}
\newcommand{\bfu}{\mbox{\boldmath$u$}}
\newcommand{\bfv}{\mbox{\boldmath$v$}}
\newcommand{\bfw}{\mbox{\boldmath$w$}}
\newcommand{\bfx}{\mbox{\boldmath$x$}}
\newcommand{\bfy}{\mbox{\boldmath$y$}}
\newcommand{\bfz}{\mbox{\boldmath$z$}}

\newcommand{\bfA}{\mbox{\boldmath$A$}}
\newcommand{\bfB}{\mbox{\boldmath$B$}}
\newcommand{\bfC}{\mbox{\boldmath$C$}}
\newcommand{\bfD}{\mbox{\boldmath$D$}}
\newcommand{\bfE}{\mbox{\boldmath$E$}}
\newcommand{\bfF}{\mbox{\boldmath$F$}}
\newcommand{\bfG}{\mbox{\boldmath$G$}}
\newcommand{\bfH}{\mbox{\boldmath$H$}}
\newcommand{\bfI}{\mbox{\boldmath$I$}}
\newcommand{\bfJ}{\mbox{\boldmath$J$}}
\newcommand{\bfK}{\mbox{\boldmath$K$}}
\newcommand{\bfL}{\mbox{\boldmath$L$}}
\newcommand{\bfM}{\mbox{\boldmath$M$}}
\newcommand{\bfN}{\mbox{\boldmath$N$}}
\newcommand{\bfO}{\mbox{\boldmath$O$}}
\newcommand{\bfP}{\mbox{\boldmath$P$}}
\newcommand{\bfQ}{\mbox{\boldmath$Q$}}
\newcommand{\bfR}{\mbox{\boldmath$R$}}
\newcommand{\bfS}{\mbox{\boldmath$S$}}
\newcommand{\bfT}{\mbox{\boldmath$T$}}
\newcommand{\bfU}{\mbox{\boldmath$U$}}
\newcommand{\bfV}{\mbox{\boldmath$V$}}
\newcommand{\bfW}{\mbox{\boldmath$W$}}
\newcommand{\bfX}{\mbox{\boldmath$X$}}
\newcommand{\bfY}{\mbox{\boldmath$Y$}}
\newcommand{\bfZ}{\mbox{\boldmath$Z$}}

\fontsize{12pt}{24pt}

\selectfont

\title {Slow Diffeomorphisms of a Manifold with $\bbT^2$ Action}

\author{
Ronit Fuchs\thanks{This paper is the author's M.Sc. thesis, being
carried out under the supervision of Prof. Leonid Polterovich, at
Tel-Aviv university.}
\\
School of Mathematical Sciences \\
Tel Aviv University
69978 Tel Aviv, Israel \\
\tt ronit.fuchs@gmail.com}

\maketitle

\numberwithin{equation}{section}

\pagestyle{myheadings}

\markboth{Ronit Fuchs}{Slow Diffeomorphisms Of A Manifold With
$\bbT^2$ Action}

\begin{abstract}
The uniform norm of the differential of the $n$-th iteration of a
diffeomorphism is called the growth sequence of the
diffeomorphism. In this paper we show that there is no lower
universal growth bound for volume preserving diffeomorphisms on
manifolds with an effective $\bbT^2$ action by constructing a set
of volume-preserving diffeomorphisms with arbitrarily slow growth.
\end{abstract}

\section{Introduction}
Let $M$ be a smooth compact connected manifold. Let $f$ be a
diffeomorphism of the manifold $M$. Define $\Gamma_n\colon$
Diff$(M) \to \bbR$, \emph{the growth sequence}, as
$$
\Gamma_n(f)=\max \left\{\max_{x\in M}\|d_xf^n\|\,,\,\max_{x\in
M}\|d_xf^{-n}\|\right\},\quad n\in \bbN.
$$
Here $f^n$ is the $n$-th iteration of $f$ and $\|d_xf\|$ is the
operator norm of the differential, calculated with respect to a
Rimannian metric on $M$. We write $a_n\succeq b_n$, when $a_n$ and
$b_n$ are two positive sequences, and there exists $c>0$ such that
$a_n \geq c b_n$ for all $n\in \bbN$. Two sequences $a_n$ and
$b_n$ are called equivalent if $a_n \succeq b_n$ and $b_n \succeq
a_n$. Under this definition, the equivalence class of the growth
sequence is an invariant of $f$ under conjugations in Diff$(M)$.
It is called \emph{the growth type} of $f$.

The growth type of a diffeomorphism is a basic dynamical invariant
(see~\cite{KH}). The behavior of \emph{the growth sequence} of
different categories of diffeomorphisms is an interesting topic.
This topic was first brought up by D'Ambra and Gromov~\cite{DAG}.
In this paper we show that there is no lower universal growth
bound for volume preserving diffeomorphisms of manifolds with an
effective $\bbT^2$ action. An action is called effective if the
only element of the group that defines the identity diffeomorphism
is the identity element.
\begin{theo}[Main Theorem]\label{mainthm}
Let $M$ be a smooth compact connected oriented manifold with an
effective $\bbT^2$ action, $\phi\colon\mathbb{T}^2\times M \to M$.
Let $\Psi$ be a positive, unbounded increasing, function on
$\mathbb{R}_+$ such that $\Psi (x)=o(x)$ for $x\rightarrow
\infty$. Then there exists a volume preserving diffeomorphism,
$f$, of $M$, such that
\begin{equation}\label{in_gamma}
0<\limsup_{n\rightarrow \infty} \frac{\Gamma_n(f)}{\Psi
(n)}<\infty.
\end{equation}
\end{theo}
We will refer to such an $f$ as to a slow diffeomorphism.

In works of Polterovich and Sikorav (see~\cite{P1, PSi}) it has
been found that there are lower growth bounds for Hamiltonian
diffeomorphisms of symplectic manifolds, $M$, with $\pi_2(M)=0$. A
Hamiltonian diffeomorphism always has fixed points and vanishing
flux (The definition of the flux is given in Section
\ref{Sec_Flux}).

In the case of symplectic, but non-Hamiltonian, diffeomorphisms
(i.e. with non-vanishing flux) Polterovich prove the existence of
lower growth bounds if the diffeomorphism has a fixed point with
some special property.

In recent works of Polterovich (see \cite{P2}) and Borichev (see
\cite{B}) it has been found that there are no lower growth bounds
("continuous spectrum") in the case of symplectic diffeomorphisms
without fixed points. They gave examples of sequences of
diffeomorphisms on the two dimensional torus with arbitrarily slow
growth.

In the case of smooth category Polterovich and Sodin ~\cite{PSo}
show that there are no growth bounds:
\begin{theo}
Let $\{u(n)\}$ be a sequence of positive real numbers which goes
to infinity as $n\rightarrow +\infty$. Then there exists a
diffeomorphism $f \in $ Diff${}_0(M) \verb|\| \{1\!\!1\}$ with a
fixed point so that $\liminf_{n\rightarrow \infty}
\frac{\Gamma_n(f)}{u(n)}<\infty$
\end{theo}

Let us emphasize that the diffeomorphisms constructed in this
theorem are dissipative, that is to say, they do not preserve any
smooth volume form.


In this work we tackle the open problem: what happens in
volume-preserving category?

Main theorem \ref{mainthm} shows that in this case, on manifolds
with an effective $\bbT^2$ action, there is no universal lower
growth bound for volume preserving diffeomorphisms.

Moreover, in certain situations, the slow diffeomorphisms
appearing in the main theorem have features similar to those of
Hamiltonian diffeomorphisms: they have fixed points and their flux
vanishes.

\begin{theo}\label{trm_slow_diff}
There exist volume-preserving slow diffeomorphisms with vanishing
flux and fixed points on the manifold $S^1\times S^2$.
\end{theo}

\subsection{Organization of the Work}

In Section \ref{Sec_MainThm}, we construct slow diffeomorphisms of
a manifold with an effective $\bbT^2$ action and thus prove
Theorem \ref{mainthm}. Farther, we give examples of slow
diffeomorphisms of manifold $S^1\times S^2$.

In Section \ref{Sec_Flux}, we define the flux homomorphism and
prove Theorem \ref{trm_slow_diff}.

\newpage
\section{Proof of Main Theorem}
\label{Sec_MainThm}
\subsection{Topological preliminaries}
\label{ManTop}

Let $M$ be a connected manifold of dimension $n$ with an effective
$\bbT^2$ action, $\phi\colon \bbT^2 \times M \to M$. Denote
$\phi(g,x)=\phi_g(x)=gx$. Consider the torus as the group $\bbT^2
= \bbR^2/\bbZ^2$.

We say that the action is free at $x\in M$, if the map $\bbT^2\to
M$, $g\mapsto \phi(g,x)$ is an embedding. We write $M_e$ as the
set of elements of $M$ where the action is free. Then, by
~\cite[Corollary B.48]{GGK}, when the action is effective, $M_e$
is open and dense.

Let $x_0\in M_e$, then the orbit $V =\bbT^2 x_0$ is diffeomorphic
to $\bbT^2$. Then there exists a neighborhood, $W$, of $x_0$ such
that $\phi(\cdot\:,x)\colon \bbT^2 \to M$ is an embedding for all
$x$ in the neighborhood. Let $D$ be an open disc in $W$ of
dimension $n-2$ such that $D$ is transversal to $V$. Finally,
denoted by $U$ the orbit of the disc under the torus action, then
$U = \bbT^2 D \cong \bbT^2\times D$.

Let $u=(u_1,\ldots,u_{n-2})$ where $\sum_1^{n-2} u_i^2<1$, be the
coordinates of the disc $D$. Let $(\varphi_1,\varphi_2)$ where
$\varphi_1,\varphi_2\in S^1=\bbR/\bbZ$ be coordinates of the torus.
Accordingly, the triple $(\varphi_1,\varphi_2,u)$ represents
coordinates in $U$.

Let $\tilde{B}$ and $B$ be two sets satisfying the following:
\begin{itemize}
\item $\tilde{B} \subset B \subset D$
\item $B$ is a compact set
\item $\tilde{B}$ is an open set.
\end{itemize}
We define a smooth function $A\colon D\to \bbR$ such that
$\supp(A) \subseteq B$ and $A|_{\tilde{B}} = 1$.

\subsection{Constructing the Diffeomorphism}
\label{proof}

We define the diffeomorphism on $M\setminus U$ and on $U$
separately. On the set $U$ we construct a slow diffeomorphism
using the function from Borichev's theorem~\cite{B} and the
function $A\colon D\to \bbR$.

For a function $F\colon S^1\to \bbR$ and $\alpha\in \bbR$, we
consider the Weyl sum
$$
W(N,x,\alpha) = \sum_{k=0}^{N-1}F(x + k \alpha).
$$
\begin{theo}[Borichev]
Let $\Psi$ be a positive, unbounded increasing, function on
$\mathbb{R}_+$ such that $\Psi (x)=o(x)$ for $x\rightarrow
\infty$. Then, there exists a real-analytic and 1-periodic
function, $F\colon \bbR \to \bbR$, and $\alpha\in \bbR$ such that
$\int_0^1 F(x)dx=0$ and
\begin{equation}\label{in_Weyl}
0<\limsup_{N\rightarrow \infty} \max_{0\leq x<1}
\frac{W'(N,x,\alpha)}{\Psi (N)}<\infty.
\end{equation}
\end{theo}

Let $F$ be the function and $\alpha$ be the constant from
Borichev's theorem. Then, we define $f_1\colon U \to U$ as the
following diffeomorphism:
$$
f_1(x)=(\varphi_1 + \alpha\,,\,\varphi_2 + A(u)\cdot
F(\varphi_1)\,,\,u),
$$
where we use the coordinates $(\varphi_1,\varphi_2,u)=x\in U$
previously mentioned.

On the set $M \setminus U$ we define $f_2\colon M \setminus U \to
M \setminus U$ as the action of the element $(\alpha,0)\in \bbT^2$
on $M \setminus U$: $f_2(x)=\phi_{(\alpha,0)}x$. Indeed, $U$ is
the orbit of a set, hence, $f_2$ is onto $M\setminus U$.

Denote $f\colon M \to M$ as 
$$
f(x)=\left \{
    \begin{array}{c l}
    f_1(x)      & x\in U \\
    f_2(x)      & x\in M\setminus U. \\
  \end{array}
\right.
$$
It is clear that $f$ is a diffeomorphism. We show that $f$
satisfies the conditions in the main theorem(\ref{mainthm}).

\begin{lem}\label{lem_U}
The diffeomorphism $f_1$ satisfies the inequality
$$
0<\limsup_{n\rightarrow \infty} \frac{\Gamma_n(f)}{\Psi
(n)}<\infty.
$$
on the submanifold $U$.
\end{lem}

\begin{proof}
Take $(\varphi_1,\varphi_2,u)$ as the coordinates on $U$. The
$m$-th iteration of $f_1$ is equal to
$$
f_1^m=(\varphi_1+m\alpha\,,\,\varphi_2 + A(u)\cdot
W(m,\varphi_1,\alpha)\,,\,u).
$$
Hence,
$$
d_xf_1^m = \left(
  \begin{array}{c c c c c}
    1       & 0      & 0 & \cdots & 0 \\
    A(u)\cdot W'(m,\varphi_1,\alpha)
            & 1      &\frac{\partial A}{\partial u_1}\cdot W(m,\varphi_1,\alpha)
                            & \cdots & \frac{\partial A}{\partial u_{n-2}}\cdot W(m,\varphi_1,\alpha)  \\
    0       & 0      & 1    &        &   \\
    \vdots  & \vdots &      & \ddots &   \\
    0       & 0      &      &        & 1 \\
  \end{array}
  \right)
$$
where $u=(u_1,\cdots,u_{n-2})\in D$.

Define a norm of a $n\times n$ matrix $Q=(q_{ij})$ as
\begin{eqnarray}\label{norm}
\|Q\|=\max_i{\sum_j|q_{ij}|}.
\end{eqnarray}
Every other norm is equivalent to
$\|\cdot\|$, hence, by using the fact that $A(u)$ and its
derivatives are bounded functions of $u$, it is sufficient to
prove that $F$ and $\alpha$ satisfy the following condition:
$$
0<\limsup_{N\rightarrow
\infty}\max_{0\leq\varphi<1}\frac{\left|W(N,\varphi,\alpha)\right|+\left|W'(N,\varphi,\alpha)\right|}{\Psi(N)}<\infty.
$$

The lower bound directly follows inequality (\ref{in_Weyl}). For
the upper bound we use the fact that there exists $\tilde{x}\in
[0,1)$ such that $W(N,\tilde{x},\alpha) = 0$. Indeed, $\int_0^1
W(N,x,\alpha)dx=0$ and $W(N,x,\alpha)$ is continuous as a function of
$x$. Hence, using $W(N,\varphi,\alpha)=\int_{\tilde{x}}^\varphi
W'(N,x,\alpha)dx + W(N,\tilde{x},\alpha)$, we get
\begin{eqnarray*}
\max_{0\leq\varphi<1}
\left|W(N,\varphi,\alpha)\right|&=&\max_{0\leq\varphi<1}
\left|\int_{\tilde{x}}^\varphi
W'(N,x,\alpha)dx\right| \leq {} \\
{}&\leq& \max_{0\leq\varphi<1} \left|W'(N,\varphi,\alpha)\right|.
\end{eqnarray*}
Hence, inequality (\ref{in_Weyl}) yields
$$
\limsup_{N\rightarrow
\infty}\max_{0\leq\varphi<1}\frac{\left|W(N,\varphi,\alpha)\right|}{\Psi(N)}<\limsup_{N\rightarrow
\infty}\max_{0\leq\varphi<1}\frac{\left|W'(N,\varphi,\alpha)\right|}{\Psi(N)}<\infty
$$
and we prove the upper bound in $U$.
\end{proof}

\begin{lem}\label{lem_M-U}
Look at $\Gamma_m(f_2)$ as a function of $m$. Then,
$\Gamma_m(f_2)$ is bounded.
\end{lem}

\begin{proof}
The diffeomorphism $f_2 \colon M\setminus U \to M\setminus U$ is
the action of $(\alpha,0)\in \bbT^2$ on $M\setminus U$ and the two
dimensional torus is compact, therefore, $\|d_xf_2^m\|$ is
bounded.
\end{proof}

%

Using Lemmas \ref{lem_U} and \ref{lem_M-U} we get that the
diffeomorphism $f$ satisfies the inequality
$$
0<\limsup_{n\rightarrow \infty} \frac{\Gamma_n(f)}{\Psi
(n)}<\infty.
$$

The manifold $M$ is oriented, hence, there exists a volume form on
$M$. Define the volume form $\omega$ on $M$ as follows. Let
$\Omega_0=d\varphi_1 \wedge d\varphi_2 \wedge du_1 \wedge \ldots
\wedge du_{n-2}$ on $U$. Let $\Omega$ be any extension of
$\Omega_0$ to the entire manifold. Let $\mu$ be a Haar measure on
the torus. We average the pullbacks of $\Omega$ by the
diffeomorphism $\phi_g$, where $g \in \mathbb{T}^2$ and set
$$
\omega = \int_{\bbT^2} g^* \Omega d\mu(g).
$$

\begin{lem}
The diffeomorphism $f$ preserves the volume form $\omega$.
\end{lem}

\begin{proof}
The action preserves $\omega$, for every $g \in \mathbb{T}^2$
$g^*\omega = \omega$. The diffeomorphism $f_2$ is the action of
the element $(\alpha,0)\in \bbT^2$ on $M$, hence, $f$ preserves
$\omega$ on $M\setminus U$.

On $U$ $\det (d_xf_1)=1$. Hence $f$ preserves the volume form
$\omega|_U=\Omega_0=d\varphi_1 \wedge d\varphi_2 \wedge du_1
\wedge \ldots \wedge du_{n-2}$ on $U$.  Therefore, $f$ preserves
the volume form $\omega$ on $M$ as required.
\end{proof}

\subsection{Examples of Manifolds with Slow Diffeomorphisms}
\label{Example}

An example of manifolds satisfying the conditions of the main
theorem are the spheres, $S^n$, where $n\geq 3$.

Another example is $S^1\times S^2$. We will construct two
diffeomorphisms that satisfies inequality (\ref{in_gamma}) on
$S^1\times S^2$.



Consider the unit sphere $S^2=\{x^2+y^2+z^2=1\}\subset \bbR^3$.
Define the polar coordinates $(\rho,\Theta)$ in the $(x,y)$-plane:
$x=\rho cos \Theta$ and $y=\rho sin \Theta$. Put
$\theta=\frac{\Theta}{2\pi}$. Then $\{(\theta,z): \theta\in [0,1),
z\in (-1,1)\}$ can be taken as coordinates of the sphere without
the poles.

Let $R_\alpha$ be the rotation of $S^2$ around the $z$-axis by
angle $2\pi\alpha \in S^1=\bbR/2\pi \bbZ$:
$$R_\alpha(w)=(\theta+\alpha,z)$$
where $w=(\theta,z)$.

\begin{examp} \label{Example_1}
\rm{Consider the torus $\bbT^2=\bbR^2/\bbZ^2$ and define a torus
action on $S^1\times S^2$ as follows
$$\phi_{(\varphi_1,\varphi_2)}(\lambda,w)=(\lambda+\varphi_1,R_{\varphi_2}(w))$$
where $(\varphi_1,\varphi_2)\in \bbT^2$, $w\in S^2$ and $\lambda
\in S^1=\bbR/\bbZ$.

Let us construct the set $U$ from Section \ref{ManTop}. The action
$\phi$ is free on the element $(0,w_0)\in S^1\times S^2$, where
$w_0=(0,0)\in S^2$. The set $D=\{(0,w): w=(0,z),
z\in(-1,1)\}\subset S^1 \times S^2$ is a one dimensional disc.
Then $U=\bbT^2D$ and we get $U=S^1\times (S^2\backslash
\{a_1,a_2\})$, where $a_1$ and $a_2$ are the poles of the sphere.

Let $\Psi$ be a positive, unbounded increasing, function on
$\mathbb{R}_+$ such that $\Psi (x)=o(x)$ for $x\rightarrow
\infty$. Let $F$ be the function from Borichev's theorem and
$\alpha$ be the constant from Borichev's theorem. In this case the
slow diffeomorphism of $S^1\times S^2$, satisfying inequality
(\ref{in_gamma}), is
$$
f(\lambda,w)=\left \{
    \begin{array}{l l}
    f_1(\lambda,w)=(\lambda+\alpha,\,R_{A(z)F(\lambda)}(w)) & (\lambda,w)\in U,\quad w=(\theta,z) \\
    f_2(\lambda,w)=(\lambda+\alpha,w)                       & (\lambda,w)\in S^1 \times \{a_1,a_2\}. \\
  \end{array}
\right.
$$
Notice the diffeomorphism does not have a fixed point.}
\end{examp}

\begin{examp} \label{Example_2}
\rm {Define the torus action $\widetilde{\phi} \colon \bbT^2
\times (S^1\times S^2) \to S^1\times S^2$ as follows
$$\widetilde{\phi}_{(\varphi_1,\varphi_2)}(\lambda,w)=(\lambda+\varphi_2,R_{\varphi_1}(w))$$
where $(\varphi_1,\varphi_2)\in \bbT^2$, $w\in S^2$ and $\lambda
\in S^1$.

Define the set $U$ as in Example \ref{Example_1}. The slow
diffeomorphism in this case is
$$
\widetilde{f}(\lambda,w)=\left \{
    \begin{array}{l l}
    \widetilde{f_1}(\lambda,w)=(\lambda+A(z)F(\theta),\,R_{\alpha}(w)) & (\lambda,w)\in U,\quad w=(\theta,z)\\
    \widetilde{f_2}(\lambda,w)=(\lambda,w)                             & (\lambda,w)\in S^1 \times \{a_1,a_2\} \\
  \end{array}
\right.
$$
Notice $\widetilde{f}$ has fixed points at $S^1 \times
\{a_1,a_2\}$.}
\end{examp}

\newpage
\section{Growth and Flux}
\label{Sec_Flux}

Let $(M,\omega)$ be a closed manifold of dimension $n$ with a
volume form, $\omega$. Let $\textrm{\emph{Diff}}_0(M,\omega)$ be
the group of volume preserving diffeomorphisms isotopic to the
identity. First, let us define the flux homomorphism on
$\pi_1(\textrm{\emph{Diff}}_0(M,\omega))$.

For any $(n-1)$-cycle $C$ and any $\{f_t\}\in
\pi_1(\textrm{\emph{Diff}}_0(M,\omega))$ we define an $n$-cycle
$\Phi(C) = \cup_t f_t(C)$ in $M$.

Define the \emph{flux} homomorphism
$$
\overline{flux_\omega} \colon
\pi_1(\textrm{\emph{Diff}}_0(M,\omega)) \to H^{n-1}(M,\bbR)
$$
as follows. Choose a loop $\{f_t\}$ representing element
$\alpha\in \pi_1 (Diff_0(M,\omega))$. Put
$$
(\overline{flux_\omega}([\{f_t\}]),[C])= (\omega,\Phi(C))
\qquad\textrm{ for all } [C]\in H_{n-1}(M,\bbR).
$$
One can show that this definition does not depend on choice of a
loop $\{f_t\}$ representing $\alpha$ and the cycle $C$
representing the homology class.

The image $\Gamma =
\overline{flux_\omega}(\pi_1(\textrm{\emph{Diff}}_0(M,\omega)))
 \subset H^{n-1}(M,\bbR)$ is called the \emph{flux group}.

The notion of \emph{flux} can be extended to diffeomorphisms $f\in
\textrm{Diff\,}_0(M,\omega)$. Define
$$
flux_\omega \colon \textrm{\emph{Diff}}_0(M,\omega) \to
H^{n-1}(M,\bbR)/\Gamma
$$
as follows. Choose any path $\{f_t\}$ of volume preserving
diffeomorphisms with $f_0=id$ and $f_1=f$. Put $\Phi(C) = \cup_t
f_t(C)$. Then $flux_\omega$ is given by
$$
(flux_\omega(f),[C])= (\omega,\Phi(C)) \qquad\textrm{ for all }
[C]\in H_{n-1}(M,\bbR).
$$

Using the fact that $f_t$ is volume preserving for every $t$,
$(\omega,\Phi(C))$ does not depend on the choice of the element
$C\in [C]$.

For $f \in \textrm{\emph{Diff}}_0(M,\omega)$ we can choose
different paths, $\{f_t\}$ and $\{g_t\}$, with $f_0=g_0=id$ and
$f_1=g_1=f$. However, the difference between the fluxes of these
paths connecting $1\!\!1$ to $f$ belongs to $\Gamma$ and thus the
flux is well defined.

Let us return to the example of $(S^1\times S^2,\omega)$, when
$\omega$ is the volume form constructed in Section \ref{proof},
after normalization, $\int_{S^1 \times S^2}\omega = 1$. We
calculate the flux of the slow diffeomorphisms $f$ and
$\widetilde{f}$. First, let us calculate the \emph{flux group},
$\Gamma$, of this manifold. Make the following identifications:
$$
H_2(S^1\times S^2,\bbZ)=\bbZ,\quad H_3(S^1\times S^2,\bbZ)=\bbZ,
\quad  H^2(S^1\times S^2,\bbZ)=\bbZ\subset\bbR=H^2(S^1\times
S^2,\bbR).
$$
For each $\gamma\in \pi_1(\textrm{\emph{Diff}}_0(M,\omega))$ and
$[a]\in H_2(S^1\times S^2,\bbZ)$ we have
$$
(flux_\omega(\gamma),a)=(\omega,\Phi(a)).
$$
Here $\Phi$ is a functional from $\bbZ$ to $\bbZ$ and the value of
$\omega$ on the generator of $H_3(S^1\times S^2,\bbZ)$ equals $1$.
Hence, $(\omega,\Phi a)\in \bbZ$ and $flux(\gamma)\in
H^2(S^1\times S^2,\bbZ)$ which implies that $\Gamma\subset \bbZ$.
On the other hand, let us look at the following loop $\gamma\in
\pi_1(\textrm{\emph{Diff}}_0(M,\omega))$ :
$$
\gamma(\lambda,w)(t)=(\lambda+t,w)
$$
where $\lambda\in S^1$, $w\in S^2$ and $t\in [0,1]$. Fix
$\lambda_0\in S^1$. Let $[C]=[\lambda_0\times S^2]$ be the
generator of $H_2(S^1\times S^2,\bbZ)$. Then $\Phi(C)=S^1\times
S^2$. Therefore,
$$
(flux(\gamma),[C])=(\omega,\Phi(C))=1,
$$
and we conclude that $\bbZ \subset \Gamma$. Hence $\Gamma=\bbZ$.

Let us calculate the flux of the diffeomorphism from Example
\ref{Example_1}. The set $U=S^1\times (S^2\backslash \{a_1,a_2\})$
is the set constructed in Example \ref{Example_1}. As before $C$
stands for $\lambda_0\times S^2$. Then
$(\omega,\Phi(C))=(\omega,\Phi(C'))$ where $C'=C\cap U =
C\backslash \lambda_0\times \{a_1,a_2\}$, and $a_1$, $a_2$ are the
poles of $S^2$. Indeed, the dimension of $\Phi(\lambda_0\times
\{a_1,a_2\})$ is one, hence $(\omega,\Phi(\lambda_0\times
\{a_1,a_2\}))=0$.

Let $f_t(\lambda,w)=(\lambda+t\alpha,\,R_{t A(z)F(\lambda)}(w))$,
$t\in [0,1]$, where $w=(\theta,z)\in S^2$, be a path of volume
preserving diffeomorphisms on $U$ with $f_0=id$ and $f_1=f$. Then,
$$
f_t\big({\{\lambda_0\times
(S^2\backslash\{a_1,a_2\})\}}\big)=\{(\lambda_0+t\alpha)\times
(S^2\backslash\{a_1,a_2\})\}
$$
and
$$
\Phi \big(\lambda_0\times (S^2\backslash\{a_1,a_2\})\big)=
\big([\lambda_0,\lambda_0+\alpha] \times
(S^2\backslash\{a_1,a_2\})\big).
$$
Thus,
$$
(flux(\{f_t\}),[C])=(\omega,\Phi(C'))=\alpha \quad\quad(\mod 1).
$$
Hence, $flux(\{f\})([C])=\alpha \,\, (\mod 1)$.

Now, let us calculate the flux of $\widetilde{f}$ from Example
\ref{Example_2}. Let $\widetilde{f_t}(\lambda,w)=(\lambda+t
A(z)F(\theta),\,R_{t \alpha}(w))$, $t\in [0,1]$, where
$w=(\theta,z)\in S^2$, be a path of volume preserving
diffeomorphisms with $\widetilde{f_0}=id$ and
$\widetilde{f_1}=\widetilde{f}$. Then,
$$
\widetilde{f_t}\big({\{\lambda_0\times
(S^2\backslash\{a_1,a_2\})\}}\big)= \{(\lambda,\theta,z) \in S^1
\times S^1 \times (-1;1)\,|\, \lambda = \lambda_0+t A(z)
F(\theta)\}
$$
and $\Phi\big(\lambda_0\times (S^2\backslash\{a_1,a_2\})\big)$ is
a domain bounded by hypersurfaces $\{\lambda=\lambda_0\}$ and
$\{\lambda=\lambda_0+A(z)F(\theta)\}$.

Thus,
$$
(flux(\{\widetilde{f_t}\}),[C])=(\omega,\Phi(C'))=\int_{-1}^1
\int_0^1 A(z)F(\theta) d\theta dz \quad\quad(\mod 1)
$$
where $C=\lambda_0\times S^2$.

Now, from Borichev's theorem, $\int_0^1 F(\theta) d\theta = 0$.
Hence, $flux(\{\widetilde{f}\})([C])\equiv 0 \,\,(\mod 1)$.

Since $H_1(S^2\times S^1,\bbZ)=\bbZ$ we get that $flux(\{f_t\}) =
0 \,\,(\mod 1)$.

In conclusion, in Example \ref{Example_2} we have a diffeomorphism
with fixed points and vanishing flux. That proves Theorem
\ref{trm_slow_diff}.

\end{document}